\def\E{\end{document}}
\documentclass[11pt]{article}
\usepackage{amssymb,amsmath}

\begin{document}
\title{Some Results on the Scattering Theory for a
Schr\"{o}dinger Equation with Combined Power-Type
Nonlinearities}
 \author{Xianfa Song {\thanks{E-mail: songxianfa2004@163.com(or
 songxianfa2008@sina.com)
}}\\
\small Department of Mathematics, School of Science, Tianjin University,\\
\small Tianjin, 300072, P. R. China }

\maketitle
\date{}

\newtheorem{theorem}{Theorem}[section]
\newtheorem{definition}{Definition}[section]
\newtheorem{lemma}{Lemma}[section]
\newtheorem{proposition}{Proposition}[section]
\newtheorem{corollary}{Corollary}[section]
\newtheorem{remark}{Remark}[section]
\renewcommand{\theequation}{\thesection.\arabic{equation}}
\catcode`@=11 \@addtoreset{equation}{section} \catcode`@=12

\begin{abstract}

In this paper, we consider the Cauchy problem
\begin{align*} \left\{\begin{array}{ll}&i
u_t+\Delta u=\lambda_1|u|^{p_1}u+\lambda_2|u|^{p_2}u, \quad t\in\mathbb{R}, \quad x\in\mathbb{R}^N\\
&u(0,x)=\varphi(x)\in \Sigma, \quad x\in\mathbb{R}^N,
\end{array}
\right. \end{align*} where $N\geq 3$,
$0<p_1<p_2\leq\frac{4}{N-2}$, $\lambda_1\in\mathbb{R}\setminus\{0\}$ and $\lambda_2\in\mathbb{R}$ are
constants, $\Sigma=\{f\in
H^1(\mathbb{R}^N); |x|f\in L^2(\mathbb{R}^N)\}$. Using the strategy in \cite{Cazenave2, Cazenave3} and  taking some elementary techniques which differ from the pseudoconformal conservation law, we obtain some scattering properties, which partly solve the open problems of Terence Tao, Monica Visan and Xiaoyi Zhang[The nonlinear
Schr\"{o}dinger equation with combined power-type nonlinearities,
Communications in Partial Differential Equations, 32(2007), 1281--1343].
As a byproduct, we establish the scattering theory in $\Sigma$ for
\begin{align*} \left\{\begin{array}{ll}&i
u_t+\Delta u=\lambda|u|^pu, \quad t\in\mathbb{R}, \quad x\in\mathbb{R}^N\\
&u(0,x)=\varphi(x), \quad x\in\mathbb{R}^N
\end{array}
\right. \end{align*}
with $\lambda>0$ and $\frac{2}{N}<p<\alpha_0$ with $\alpha_0=\frac{2-N+\sqrt{N^2+12N+4}}{2N}$, which is also an open problem in this direction.



\end{abstract}

\section{Introduction}
\qquad In this paper, we consider the following
Cauchy problem
\begin{align}
\label{system0} \left\{\begin{array}{ll}&i u_t+\Delta
u=\lambda_1|u|^{p_1}u+\lambda_2|u|^{p_2}u, \quad t\in\mathbb{R}, \quad x\in\mathbb{R}^N\cr
 &u(0,x)=\varphi(x), \quad x\in\mathbb{R}^N,
\end{array}
\right. \end{align} where $N\geq 3$, $0<p_1<p_2\leq \frac{4}{N-2}$,
$\lambda_1\in\mathbb{R}\setminus\{0\}$ and $\lambda_2\in\mathbb{R}$ are constants. The model (\ref{system0}) appears in the theory of Bose-Einstein
condensation and nonlinear optics(see \cite{Cazenave2,Kato, Strauss}. We are interested in the scattering
properties of the solutions. In convenience, we take
the same conventional notions of scattering theory as those in \cite{Cazenave2} below.

Let $I$ be an interval containing $0$, Duhamel's formula implies
that $u$ is a solution of (\ref{system0}) on $I$ if and only if
$u$ satisfies
\begin{align}
u(t)=\mathcal{J}(t)\varphi-i\int_0^t\mathcal{J}(t-s)\lambda_1|u(s)|^{p_1}u(s)ds
-i\int_0^t\mathcal{J}(t-s)\lambda_2|u(s)|^{p_2}u(s)ds\label{101231}
\end{align}
for all $t\in I$, where $\mathcal{J}(t)=e^{it\Delta}$ is the one
parameter group generated by the free Schr\"{o}dinger equation.
Let $X$ be a Banach space -- $X$ can be $\Sigma$,
$H^1(\mathbb{R}^N)$ or $L^2(\mathbb{R}^N)$ in this paper. Here
the pseudoconformal space
\begin{align} \Sigma:=\{f\in
H^1(\mathbb{R}^N); |x|f\in L^2(\mathbb{R}^N)\}\quad {\rm with \
norm}\ \|f\|_{\Sigma}=\|f\|_{H^1_x}+\|xf\|_{L^2_x}.
\label{10243}\end{align} Assume that the solution
$u_{\varphi}(t,x)$ is defined for all $t\geq 0$ with initial value
$\varphi\in X$. We say that $u_+$ is the scattering state of
$\varphi$ at $+\infty$ if the limit
\begin{align}
u_+=\lim_{t\rightarrow +\infty}
\mathcal{J}(-t)u_{\varphi}(t)\label{10241}
\end{align}
exists in $X$. Similarly, we say that $u_-$ is the scattering
state of $\varphi$ at $-\infty$ if the limit
\begin{align}
u_-=\lim_{t\rightarrow -\infty}
\mathcal{J}(-t)u_{\varphi}(t)\label{10242}
\end{align}
exists in $X$.

About the topic of
scattering theory, there are many results on the Cauchy problem
\begin{align}
\label{classic1} \left\{\begin{array}{ll}&i u_t+\Delta
u=\lambda|u|^pu, \quad t\in \mathbb{R}, \quad x\in\mathbb{R}^N\cr
 &u(0,x)=\varphi(x), \quad x\in\mathbb{R}^N.
\end{array}
\right. \end{align}  Different scattering theories had been
constructed in many papers. First, we review the results on (\ref{classic1}) with $0<p<\frac{4}{N-2}$.
If $\lambda>0$ and $p\leq
\frac{2}{N}$, then there are no nontrivial solution of (\ref{classic1})
has scattering states, even for $L^2(\mathbb{R}^N)$ topology(see\cite{Barab, Strauss, Tsutsumi1}).
If $\lambda>0$, $p>\frac{2}{N}$ and $\varphi\in H^1(\mathbb{R}^N)$, scattering theory in the energy space $H^1(\mathbb{R}^N)$ was established(see \cite{Ginibre1, Ginibre2, Ginibre3, Nakanishi1, Nakanishi2}). If $\lambda>0$, $p>\frac{4}{N+2}$ and $\varphi\in \Sigma$, a low energy scattering theory exists
in $\Sigma$, especially, if $p>\alpha_0$ with $\alpha_0=\frac{2-N+\sqrt{N^2+12N+4}}{2N}$, scattering theory always exists in $\Sigma$(see \cite{Ginibre4, Nakanishi3, Tsutsumi1}). However, we don't know whether $u_{\pm}\in \Sigma$ if $\lambda>0$
and $\frac{2}{N}<p<\alpha_0$ with $\|\varphi\|_{\Sigma}$ is large, which is also an open problem in this direction.  For
the case of $\lambda<0$, there is no low energy scattering if
$p<\frac{4}{N+2}$. If $\lambda<0$ and $\frac{4}{N+2}<p<\frac{4}{N}$, then a low
energy scattering theory exists in $\Sigma$.  If
$\lambda<0$ and $p\geq \frac{4}{N}$, then some solutions will blow up
in finite time, some solutions with small initial data in
$H^1(\mathbb{R}^N)$ are global and bounded in $H^1(\mathbb{R}^N)$
(see \cite{Fang, Cazenave2, Cazenave3,
 Ginibre1, Glassey, Lin, Strauss, Tsutsumi1, Weinstein} and the
references therein). Very recently, there are many results on the scattering for (\ref{classic1}) with  $p=\frac{4}{N-2}$(the energy-critical case). In \cite{Tao}, Tao dealt with global well-posedness and scattering for the higher-dimensional energy-critical nonlinear Schr\"odinger equation for radial data. For the focusing case($\lambda<0$), we can refer \cite{KM1,KM2, Killip,Li} to see the results on global well-posedness, scattering and blow-up. For the defocusing case($\lambda>0$), we can refer \cite{Colliander,Ryckman,Visan} to see more information on the topic.

The immediate motion of this paper is \cite{Tao}. Recently, T. Tao, M. Visan and X. Y. Zhang studied the scattering
properties of (\ref{system0}) with large initial data in the energy space
$H^1(\mathbb{R}^N)$ and in $\Sigma$. Their results were summarized in Table 1 of \cite{Tao}:
\begin{tabular}{|l|}
 $\lambda_2$,\quad  $\lambda_1$ \qquad \qquad \qquad $p_1$, \ $p_2$\qquad \qquad \qquad GWP \qquad Scattering\qquad \quad Provided \\
\\
$\lambda_2>0$,\ $\lambda_1\in \mathbb{R}$ \quad $0<p_1<p_2\leq \frac{4}{N-2}$ \qquad\quad
 \checkmark \qquad\qquad ?\qquad \qquad \qquad \qquad -\\
 $\lambda_2>0,\ \lambda_1>0 \quad \frac{4}{N}\leq p_1<p_2\leq \frac{4}{N-2}$ \qquad\quad
 \checkmark \qquad\qquad i n \ $H^1_x$\qquad \qquad \qquad-\\
  $\lambda_2>0,\ \lambda_1\in \mathbb{R} \quad \frac{4}{N}\leq p_1<p_2\leq \frac{4}{N-2}$ \qquad\quad
 \checkmark \qquad\qquad i n \ $H^1_x\qquad \qquad M(u_0)<<1$\\
  $\lambda_2>0,\ \lambda_1>0 \quad \alpha(N)<p_1<p_2\leq \frac{4}{N-2}$ \qquad
 \checkmark \qquad\qquad i n \ $\Sigma$ \qquad \qquad $u_0\in \Sigma$\\
  $\lambda_2<0,\ \lambda_1\in \mathbb{R} \quad 0<p_1<p_2\leq \frac{4}{N}$ \qquad\qquad\quad
 \checkmark \qquad\qquad ?\qquad \qquad \qquad \qquad-\\
  $\lambda_2<0,\ \lambda_1>0 \quad 0<p_1<p_2, \frac{4}{N}<p_2\leq \frac{4}{N-2}\quad
 \times \qquad\quad \times \qquad y_0>0,\ E(u_0)<0$\\
 $\lambda_2<0,\ \lambda_1<0 \quad \frac{4}{N}<p_1<p_2\leq \frac{4}{N-2}\qquad\qquad
 \times \qquad\quad \times \qquad  y_0>0, \ E(u_0)<0$\\
 $\lambda_2<0,\ \lambda_1<0 \quad 0<p_1\leq\frac{4}{N}<p_2\leq \frac{4}{N-2}\qquad
 \times \qquad\quad \times \quad y_0>0,E(u_0)<-CM(u_0)$\\
\end{tabular}
From the summary above, we see that there are some open problems on the scattering properties for (\ref{system0})
in the following cases:

Case (i) $\lambda_2>0$, $\lambda_1\in \mathbb{R}$, $0<p_1<p_2\leq \frac{4}{N-2}$.

Case (ii) $\lambda_2<0$, $\lambda_1\in\mathbb{R}$,
$0<p_1<p_2<\frac{4}{N}$.

Our aim is to give some results on the scattering theory of (\ref{system0})
in the two cases above. To do this, we need some observations. If one of $\lambda_1$ and $\lambda_2$ is positive and another is negative, then one of the nonlinearities is defocusing and another is focusing, hence we need to analyze the interaction between the nonlinearity $\lambda_1|u|^{p_1}u$ and $\lambda_2|u|^{p_2}u$. Under some suitable assumptions, we obtain some new scattering properties for (\ref{system0}) and partly solve the open problems in \cite{Tao}. As a byproduct, we establish a scattering theory in $\Sigma$  for (\ref{classic1})
with $\lambda>0$ and $\frac{2}{N}<p<\alpha_0=\frac{2-N+\sqrt{N^2+12N+4}}{2N}$, which also solves an open problem in this direction. First, we establish the no scattering results as follows.

{\bf Theorem 1.} (No Scattering Results) {\it Assume that $u(t,x)$
is the nontrivial solution of (\ref{system0}) with initial value $\varphi\in \Sigma$. Then
$\mathcal{J}(-t)u(t)$ does not have any strong limit in
$L^2(\mathbb{R}^N)$ if

(i) $\lambda_1\in \mathbb{R}\setminus\{0\}$, $\lambda_2\in \mathbb{R}$, $0<p_1\leq
\frac{2}{N}$, $0<
p_2<\frac{4}{N}$ with $N\geq 3$ or

(ii) $\lambda_1\in \mathbb{R}\setminus\{0\}$, $\lambda_2\geq 0$, $0<p_1\leq
\frac{2}{N}$, $\frac{4}{N}\leq p_2\leq\frac{4}{N-2}$ with $N\geq 6$. }

The second theorem is about the scattering in $\Sigma$ for (\ref{system0}).

{\bf Theorem 2.} (Scattering in $\Sigma$) {\it Assume that $u(t,x)$
is the nontrivial solution of (\ref{system0}) with initial value
$\varphi\in \Sigma$. Then there exist $u_{\pm}\in
\Sigma$ such that
\begin{align}
\mathcal{J}(-t)u(t)\rightarrow u_{\pm} \quad {\rm in } \
\Sigma \quad {\rm as } \ t\rightarrow
\pm\infty\label{912231}
 \end{align}
 if \quad (1) $\lambda_1<0$, $\lambda_2>0$ and
$\frac{4}{N+2}<p_1<p_2<\frac{4}{N-2}$ or

(2) $\lambda_1<0$, $\lambda_2<0$ and
$\frac{4}{N+2}<p_1<p_2<\frac{4}{N}$ or

 (3) $\lambda_1>0$, $\lambda_2\geq 0$, $\frac{2}{N}<p_1<p_2\leq \frac{4}{N-2}$ or

 (4) $\lambda_1>0$, $\lambda_2<0$, $\frac{2}{N}<p_1<p_2<\frac{4}{N}$ and
\begin{align}
\|\varphi(x)\|_{L^2}^{\frac{4}{N}}
<\frac{(4-Np_1)}{2N(p_2-p_1)C_N}\left(\frac{p_2+2}{|\lambda_2|}\right)^{\frac{4-Np_1}{N(p_2-p_1)}}
\left(\frac{\lambda_1(4-Np_1)(Np_1-2)}{2(4-Np_2)(p_1+2)}\right)^{\frac{4-Np_2}{N(p_2-p_1)}}.\label{104104}
\end{align}
Here $$
C_N=\frac{N+2}{N\|W^*\|_2^{\frac{4}{N}}}
$$
and $W^*$ is the ground state solution of
$$
\Delta W-\frac{2}{N}W+|W|^{\frac{4}{N}}W=0.
$$}

As a direct consequence of Theorem 2, we get the following corollary, which solves an open problem in this direction.

{\bf Corollary 1.1.} {\it Assume that $u(t,x)$
is the nontrivial solution of (\ref{classic1}) with
 $\lambda>0$, $\frac{2}{N}<p<\alpha_0=\frac{2-N+\sqrt{N^2+12N+4}}{2N}$ and $\varphi(x)\in \Sigma$.
Then there exist $u_{\pm}\in
\Sigma$ such that
\begin{align}
\mathcal{J}(-t)u(t)\rightarrow u_{\pm} \quad {\rm in } \
\Sigma \quad {\rm as } \ t\rightarrow
\pm\infty.
 \end{align}
 }

As a special case of Theorem 2, we can obtain the following corollary

{\bf Corollary 1.2.} {\it Assume that $u(t,x)$
is the nontrivial solution of
\begin{align}
\label{special} \left\{\begin{array}{ll}&i u_t+\Delta
u=\lambda_1|u|^{p_1}u+\lambda_2|u|^{\frac{4}{N-2}}u, \quad t\in\mathbb{R}, \quad x\in\mathbb{R}^N\cr
 &u(0,x)=\varphi(x)\in \Sigma, \quad x\in\mathbb{R}^N,
\end{array}
\right. \end{align}
with $\lambda_1>0$, $\lambda_2\geq 0$ and
$\frac{2}{N}<p_1<\frac{4}{N-2}$. Then there exist $u_{\pm}\in
\Sigma$ such that
\begin{align*}
\mathcal{J}(-t)u(t)\rightarrow u_{\pm} \quad {\rm in } \
\Sigma \quad {\rm as } \ t\rightarrow
\pm\infty.
 \end{align*}
  }
By the way, the model (\ref{special}) with $\lambda_1, \lambda_2\in \mathbb{R}$ and
$N=3$ had been studied by Zhang in \cite{Zhang}. The global well-posedness, scattering for (\ref{special}) with $\lambda_2>0$ and the blowup phenomenon for (\ref{special}) with $\lambda_2<0$ were studied here.

In the course of the proof of Theorem 2, we obtain some asymptotic behavior for the solution to (\ref{system0}), which can be stated as follows.

{\bf Theorem 3.} (Asymptotic behavior) {\it Assume that $u(t,x)$
is the nontrivial solution of (\ref{system0}) with initial value
$\varphi\in \Sigma$. Then for every $2\leq r\leq \frac{2N}{N-2}$, we have
\begin{align}
\|u(t)\|_{L^r}\leq C|1+t|^{-\frac{N(r-2)}{2r}} {\quad for \ all}\quad t\in \mathbb{R}\label{113122}
\end{align}
if (i) $\lambda_1>0$, $\lambda_2\geq 0$, $\frac{2}{N}<p_1<p_2\leq \frac{4}{N-2}$ or\\
 (ii) $\lambda_1>0$, $\lambda_2<0$, $\frac{2}{N}<p_1<p_2<\frac{4}{N}$ with (\ref{104104}). }

This paper is organized as follows: In Section 2, we will give
some preliminaries. In Section 3, we will prove Theorem 1. In Section 4, we will prove Theorem 2 and Theorem 3. In the Section 5, we will give some comments on the
results of this paper.

\section{Preliminaries}
\qquad  In the sequels, we will use $C, C', C_1, c$ and so on to denote various finite positive
constants, which depend on $p_1, p_2, N, \lambda_1, \lambda_2$ and $\varphi(x)$. The exact values may vary from line to line.

Similar to Section 7.5 of \cite{Cazenave2}, we will study (\ref{system0}) by using pseudoconformal
transformation. Since we only concern the scattering properties of the solution $u(t,x)$ to (\ref{system0}), we mainly give the arguments under the assumption that the maximal existence interval of $u(t,x)$ is $[0, +\infty)$. Similarly, we can discuss the problem when the maximal existence interval of $u(t,x)$ is $(-\infty,0]$ or $(-\infty, +\infty)$.

For $(t,x)\in \mathbb{R}\times \mathbb{R}^N$, let
\begin{align}
t=\frac{s}{1-s},\ x=\frac{y}{1-s},\quad {\rm or \ equivalently}, \
s=\frac{t}{1+t},\ y=\frac{x}{1+t}.\label{912051}\end{align} For
the function $u$ defined on $(a,b)\times \mathbb{R}^N(0\leq
a<b<+\infty$ are given ), set \begin{align}
v(s,y)=(1-s)^{-\frac{N}{2}}u(\frac{s}{1-s},
\frac{y}{1-s})e^{-i\frac{|y|^2}{4(1-s)}}=(1+t)^{\frac{N}{2}}u(t,
x)e^{-i\frac{|x|^2}{4(1+t)}}\label{912052}
\end{align}
for $y\in \mathbb{R}^N$ and $\frac{a}{1+a}<s<\frac{b}{1+b}$.
Obviously, if $u$ is defined on $(0,+\infty)$, then $v$ is defined
on $(0,1)$. And $u\in C([a,b], \Sigma)$ if and only if $v\in
C([\frac{a}{1+a},\frac{b}{1+b}], \Sigma)$. And it is easy to verify the following
identities
\begin{align}
\|\nabla v(s)\|_{L^2}^2&=\frac{1}{4}\|(x+2i(1+t)\nabla )u(t)\|_{L^2}^2,\label{10371}\\
\|\nabla u(t)\|_{L^2}^2&=\frac{1}{4}\|(y-2i(1-s)\nabla )v(s)\|_{L^2}^2,\label{10372}\\
\|v(s)\|_{L^{\beta+2}}^{\beta+2}&=(1+t)^{\frac{N\beta}{2}}\|u(t)\|_{L^{\beta+2}}^{\beta+2},\quad \beta\geq 0.\label{10373}
\end{align}

After some elementary computations, we see that $u(t,x)$ satisfies (\ref{system0}) if and only if
$v(s,y)$ satisfies the Cauchy problem
\begin{align}
\label{system2} \left\{\begin{array}{ll}i v_s+\Delta_y
v&=\lambda_1(1-s)^{\frac{Np_1-4}{2}}|v|^{p_1}v+\lambda_2(1-s)^{\frac{Np_2-4}{2}}|v|^{p_2}v\cr
&:=\lambda_1h_1(s)|v|^{p_1}v+\lambda_2h_2(s)|v|^{p_2}v, \quad s>0, \quad
y\in\mathbb{R}^N\cr
 v(0,y)&=\psi(y)=\varphi(x)e^{-\frac{i|x|^2}{4}}, \quad y\in\mathbb{R}^N.
\end{array}
\right. \end{align} (\ref{system2}) equals to the following
integral equation
\begin{align}
v(s)&=\mathcal{J}(s)\psi-i\int_0^s\mathcal{J}(s-\tau)
\lambda_1h_1(\tau)|v(\tau)|^{p_1}v(\tau)d\tau\nonumber\\
&\qquad
\qquad-i\int_0^s\mathcal{J}(s-\tau)\lambda_2h_2(\tau)|v(\tau)|^{p_2}v(\tau)d\tau.\label{system3}
\end{align}

First, we need to discuss the existence of the solution to (\ref{system2}). Since $u(t,x)$ satisfies (\ref{system0}) on $(a,b)$ if and only if $v(s,y)$ solves (\ref{system2}) on $(\frac{a}{1+a}, \frac{b}{1+b})$, recalling that Proposition 3.1 and Proposition 3.2 of \cite{Tao} had established the local well-posedness for (\ref{system0}) with $0<p_1<p_2\leq \frac{4}{N-2}$, we have the following existence result.

{\bf Proposition 2.1.} {\it Assume that $0<p_1<p_2\leq \frac{4}{N-2}$ and $v(0,y)=\psi(y)=u(0,x)e^{-\frac{i|x|^2}{4}}$ with $u(0,x)\in \Sigma$ is the initial value such that
(\ref{system0}) admits a unique strong $H^1_x$-solution $u(t,x)$ defined on $[0,+\infty)$. Then (\ref{system2}) has a unique strong $H^1_y$-solution $v(s,y)$ defined on $[0,1)$.
}

Our results in this paper are based on the following observation, its proof is similar to that of Proposition 7.5.1 in \cite{Cazenave2}, we omit the details here.

{\bf Proposition 2.2.} {\it Assume that $u\in C([0,+\infty),\Sigma)$
is the solution of (\ref{system0}) and $v\in C([0,1),\Sigma)$ is defined by
(\ref{912052}), i.e., it is the corresponding solution of (\ref{system2}). Then $\mathcal{J}(-t)u(t)$ has a strong limit in $\Sigma$ as $t\rightarrow +\infty$ if and only if $v(s)$ has a
strong limit in $\Sigma$ as $s\rightarrow 1$, and in that case
\begin{align}
\lim_{t\rightarrow
+\infty}\mathcal{J}(-t)u(t)=e^{i\frac{|y|^2}{4}}\mathcal{J}(-1)v(1)\quad
{\rm in }\quad\Sigma.\label{912155}
\end{align}}

\section{The Proof of Theorem 1}
\qquad In this short section, we give the proof of  Theorem 1.

{\bf The proof of Theorem 1:} We only give the proof of it for the
case of $t\rightarrow +\infty$. The proof of the case of $t\rightarrow -\infty$
is similar. Assume that
$$\mathcal{J}(-t)u(t)\rightarrow u_+\quad {\rm in} \
L^2(\mathbb{R}^N)\quad {\rm as }\ t\rightarrow +\infty$$ by
contradiction. Consequently,
\begin{align}
\|u_+\|_{L^2}=\|u(t)\|_{L^2}=\|\varphi\|_{L^2}>0.\label{912232}
\end{align}
By the results of Proposition 2.1, we have
$$v(s)\rightarrow w\quad {\rm in} \
L^2(\mathbb{R}^N)\quad {\rm as }\ s\rightarrow 1,$$ where
$$w=\mathcal{J}(1)(e^{-i\frac{|y|^2}{4}}u_+)\neq 0.$$
Noticing that $p_1+1<p_2+1\leq 2$ under the assumptions of ours, we have
\begin{align*}
|v(s)|^{p_1}v(s)\rightarrow |w|^{p_1}w\neq 0 \quad {\rm in}\
L^{\frac{2}{p_1+1}}(\mathbb{R}^N), \\
|v(s)|^{p_2}v(s)\rightarrow |w|^{p_2}w\neq 0 \quad {\rm in}\
L^{\frac{2}{p_2+1}}(\mathbb{R}^N)
\end{align*}
as $s\rightarrow 1$. Let $\theta\in \mathcal{D}(\mathbb{R}^N)$ be
the function satisfying
\begin{align}
<i|w|^{p_1}w,\theta>=1. \label{912241}
\end{align}
Using (\ref{system2}), we have
\begin{align*}
\frac{d}{ds}<v(s),\theta>&=<i\Delta v,
\theta>+\lambda_1(1-s)^{\frac{Np_1-4}{2}}<i|v|^{p_1}v,\theta>
+\lambda_2(1-s)^{\frac{Np_2-4}{2}}<i|v|^{p_2}v,\theta>\nonumber\\
&=<iv,
\Delta\theta>+\lambda_1(1-s)^{\frac{Np_1-4}{2}}<i|v|^{p_1}v,\theta>
+\lambda_2(1-s)^{\frac{Np_2-4}{2}}<i|v|^{p_2}v,\theta>.
\end{align*}
Noticing that $v$ is bounded in $L^2(\mathbb{R}^N)$ and (\ref{912241}), we can get
\begin{align}
|\frac{d}{ds}<v(s),\theta>|&\geq
\frac{1}{2}|\lambda_1|(1-s)^{\frac{Np_1-4}{2}}-C(1-s)^{\frac{Np_2-4}{2}}-C\nonumber\\
&\geq
\frac{1}{4}|\lambda_1|(1-s)^{\frac{Np_1-4}{2}}-C\label{912242}
\end{align}
if $s$ is closed to $1$ enough. However, (\ref{912242}) implies
that $|<v(s),\theta>|\rightarrow +\infty$ as $s\rightarrow 1$
because $\frac{Np_1-4}{2}\leq -1$, which is absurd.\hfill $\Box$

\section{The Proof of Theorem 2}
\subsection{Scattering Theory in $\Sigma$ for (\ref{system0}) with $\lambda_1<0$}
\qquad In this subsection, we focus on the scattering for (\ref{system0}) with $\lambda_1<0$.

Set
\begin{align}
\mathcal{R}_+&=\{\varphi\in \Sigma: T_{\max}=+\infty \ {\rm and}\
u_+=\lim_{t\rightarrow +\infty} \mathcal{J}(-t)u_{\varphi}(t)\
{\rm
exists}\},\label{91213w1}\\
\mathcal{R}_-&=\{\varphi\in \Sigma: T_{\min}=+\infty \ {\rm and}\
u_-=\lim_{t\rightarrow -\infty} \mathcal{J}(-t)u_{\varphi}(t)\
{\rm exists}\}.\label{91213w2}
\end{align}
For $\varphi\in \mathcal{R}_{\pm}$, we define the operators
\begin{align}
U_{\pm}(\varphi)=\lim_{t\rightarrow \pm \infty}
\mathcal{J}(-t)u_{\varphi}(t),\label{91213w3}
\end{align}
where the limit holds in $\Sigma$. Set
\begin{align}
\mathcal{U}_{\pm}=U_{\pm}(\mathcal{R}_{\pm}).\label{912134}
\end{align}
If the mappings $U_{\pm}$ are injective, we can define the wave
operators
\begin{align}
\Omega_{\pm}=U^{-1}_{\pm}:\mathcal{U}_{\pm}\rightarrow
\mathcal{R}_{\pm}.\label{912135}
\end{align}
And we also introduce the sets
\begin{align}
\mathcal{O}_{\pm}=U_{\pm}(\mathcal{R}_+\cap\mathcal{R}_-).\label{912136}
\end{align}
Denote the scattering operator $\bf{S}$ by
\begin{align}
{\bf S}=U_+\Omega_-:\mathcal{O}_-\rightarrow
\mathcal{O}_+.\label{912137}
\end{align}

{\bf Proposition 4.1.} {\it Assume that $\lambda_1<0$, $\lambda_2\in \mathbb{R}$,
$\frac{4}{N+2}<p_1<p_2<\frac{4}{N-2}$. Then for every $s_0\in
\mathbb{R}$ and $\psi\in \Sigma$, there exist
$T_m(s_0,\psi)<s_0<T_M(s_0,\psi)$ and a unique maximal solution
$v\in C((T_m, T_M),\Sigma)$ of equation (\ref{system2}). And the
solution $v$ satisfies the following properties:

(i) If $T_M=1$, then
$$
\lim_{s\rightarrow 1} \inf\{\left((1-s)^{\frac{(N+2)p_1-4}{4p_1}}+
(1-s)^{\frac{(N+2)p_2-4}{4p_2}}\right)\|v(s)\|_{H^1}\}>0.
$$

(ii) $v$ depends continuously on $\psi$ in the sense of the
mapping $\psi\rightarrow T_M$ is lower semicontinuous
$\Sigma\rightarrow (0,+\infty]$ and the mapping $\psi\rightarrow
T_m$ is upper semicontinuous $\Sigma\rightarrow [-\infty,0)$. Let
$v_n$ be the solution of (\ref{system2}) with initial value
$\psi_n$. If $\psi_n\rightarrow \psi$ in $\Sigma$ as $n\rightarrow
\infty$ and if $[T_1, T_2]\in (T_m,T_M)$, then $v_n\rightarrow v$
in $C([T_1,T_2],\Sigma)$.}

{\bf Proof:} Note that the nonlinearities in (\ref{system0}) satisfy the conditions of Theorem 4.4.6 in \cite{Cazenave2}. Similar to the proofs of Theorem 4.11.1 and Theorem 4.11.2 there,
roughly, replacing $h(s)|v|^{\alpha}v$ by $h_1(s)|v|^{p_1}v+h_2(s)|v|^{p_2}v$ with  $h_1(s)=f_1(s-s_0)$ and $h_2(s)=f_2(s-s_0)$, where
\begin{align}
\label{10125x2}f_1(s)=\left\{\begin{array}{ll}&\lambda_1(1-s)^{\frac{Np_1-4}{2}},
 \quad {\rm if} \ -\infty<s<1\\
&\lambda_1, \quad {\rm if}\ s\geq 1,
\end{array}
\right. \end{align} and
\begin{align}
\label{10125x3}f_2(s)=\left\{\begin{array}{ll}&\lambda_2(1-s)^{\frac{Np_2-4}{2}}, \quad {\rm if} \ -\infty<s<1\\
&\lambda_2, \quad {\rm if}\ s\geq 1,
\end{array}
\right. \end{align}  we can get the
results of Proposition 4.1.  We omit the details here.\hfill $\Box$

Using Proposition 2.2 and Proposition 4.1, we can prove the following two theorems.

{\bf Proposition 4.2.} {\it Assume that

(i) $\lambda_1<0$, $\lambda_2>0$ and
$\frac{4}{N+2}<p_1<p_2<\frac{4}{N-2}$ or

(ii) $\lambda_1<0$, $\lambda_2<0$ and
$\frac{4}{N+2}<p_1<p_2<\frac{4}{N}$.

Then

(i) The sets $\mathcal{R}_{\pm}$ and $\mathcal{U}_{\pm}$ are open
subsets of $\Sigma$ with $0\in \mathcal{R}_{\pm}$ and
$0\in\mathcal{U}_{\pm}$.

(ii) The operators $U_{\pm}:\mathcal{R}_{\pm}\rightarrow
\mathcal{U}_{\pm}$ and
$\Omega_{\pm}:\mathcal{U}_{\pm}\rightarrow\mathcal{R}_{\pm}$ are
all bicontinuous bijections for the $\Sigma$ topology.

(iii) The sets $\mathcal{O}_{\pm}$ are open subsets of $\Sigma$
with $0\in\mathcal{O}_{\pm}$, and the scattering operator $\bf{S}$
is a bicontinuous bijection $\mathcal{O}_-\rightarrow
\mathcal{O}_+$ for the $\Sigma$ topology.}

{\bf Proof:} The proof is similar to the standard argument of Theorem 7.5.7 in \cite{Cazenave2}. We omit the details here.\hfill$\Box$

We have further results about the wave operators $\Omega_{\pm}$
which can be read as

{\bf Proposition 4.3} {\it Assume that

(i) $\lambda_1<0$, $\lambda_2>0$ and
$\frac{4}{N+2}<p_1<p_2<\frac{4}{N-2}$ or

(ii) $\lambda_1<0$, $\lambda_2<0$ and
$\frac{4}{N+2}<p_1<p_2<\frac{4}{N}$.

Then
$\mathcal{U}_{\pm}=\Sigma$. Hence the wave operators
$\Omega_{\pm}$ are bicontinuous bijections $\Sigma\rightarrow
\mathcal{R}_{\pm}$.}

{\bf Proof:} The proof is similar to the standard argument of Theorem 7.5.9 in \cite{Cazenave2}. We omit the details here.\hfill$\Box$

Now the scattering theory in $\Sigma$ for (\ref{system0}) with $\lambda_1<0$ is the direct consequence of Proposition 4.2 and Proposition 4.3.

\subsection{Scattering Theory in $\Sigma$ for (\ref{system0}) with $\lambda_1>0$}
\qquad To establish the scattering theory in $\Sigma$ for (\ref{system0}) with $\lambda_1>0$, the key step is
to deduce that $\|v(s)\|_{H^1}$ keeps bounded as $s\rightarrow 1$. The following proposition will give the estimate for $\|v(s)\|_{H^1}$.

{\bf Proposition 4.4.} {\it Assume that
$v(s,y)$ is the solution of (\ref{system2}). Then
\begin{align}
\|v(s)\|_{L^2}&\leq C\quad {\rm for \ all} \quad 0\leq s<1.\label{912264}
\end{align}
Moreover, if

(A) $\lambda_1>0$, $\lambda_2\geq 0$, $\frac{2}{N}<p_1<p_2\leq \frac{4}{N-2}$ or

(B)$\lambda_1>0$, $\lambda_2<0$, $\frac{2}{N}<p_1<p_2<\frac{4}{N}$ and (\ref{104104}),\\
then
\begin{align}
\|\nabla v(s)\|_{L^2}^2\leq C,\quad \|v(s)\|_{L^{p_1+2}}^{p_1+2}\leq C,\quad
\|v(s)\|_{L^{p_2+2}}^{p_2+2}\leq C\quad {\rm for \ all} \quad 0\leq s<1.\label{103303}
\end{align}
}

{\bf Proof:} Noticing that
$$
\frac{d}{ds}\|v(s)\|_{L^2}=0,
$$
we can obtain
$$
\|v(s)\|_{L^2}\leq C\quad {\rm for \ all} \quad 0\leq s<1.
$$

Multiplying the first equation of (\ref{system2}) by $\bar{v}_s$, integrating it on $[0,s]\times \mathbb{R}^N$ and taking the real part of the resulting expression, we
have
\begin{align}
&\frac{1}{2}\int_{\mathbb{R}^N}|\nabla v(s)|^2dy+\frac{\lambda_1(1-s)^{\frac{Np_1-4}{2}}}{p_1+2}\int_{\mathbb{R}^N}|v(s)|^{p_1+2}dy
+\frac{\lambda_2(1-s)^{\frac{Np_2-4}{2}}}{p_2+2}\int_{\mathbb{R}^N}|v(s)|^{p_2+2}dy\nonumber\\
&=\frac{1}{2}\int_{\mathbb{R}^N}|\nabla \psi|^2dy+\frac{\lambda_1}{p_1+2}\int_{\mathbb{R}^N}|\psi|^{p_1+2}dy
+\frac{\lambda_2}{p_2+2}\int_{\mathbb{R}^N}|\psi|^{p_2+2}dy\nonumber\\
&\qquad \qquad \qquad \qquad -\frac{\lambda_1(Np_1-4)}{2(p_1+2)}\int^s_0(1-\tau)^{\frac{Np_1-4}{2}-1}\int_{\mathbb{R}^N}|v(\tau)|^{p_1+2}dyd\tau\nonumber\\
&\qquad \qquad \qquad \qquad -\frac{\lambda_2(Np_2-4)}{2(p_2+2)}\int_0^s(1-\tau)^{\frac{Np_2-4}{2}-1}\int_{\mathbb{R}^N}|v(\tau)|^{p_2+2}dyd\tau.
\label{10411}
\end{align}

Case (A) $\lambda_1> 0$, $\lambda_2\geq 0$, $\frac{2}{N}<p_1<p_2\leq \frac{4}{N-2}$.
We divide it into three subcases:

Subcase (i) $\lambda_1>0$, $\lambda_2\geq 0$, $\frac{4}{N}\leq p_1<p_2\leq \frac{4}{N-2}$;

Subcase (ii) $\lambda_1>0$, $\lambda_2\geq 0$, $\frac{2}{N}<p_1<p_2\leq \frac{4}{N}$;

Subcase (iii) $\lambda_1>0$, $\lambda_2\geq 0$, $\frac{2}{N}<p_1<\frac{4}{N}\leq p_2\leq \frac{4}{N-2}$.

In subcase (i), from (\ref{10411}), we can directly obtain $$\|\nabla v(s)\|^2_2\leq C \quad {\rm for \ all }\quad  s\in [0,1).$$ Consequently, using
Gagliardo-Nirenberg's inequality and $\|v(s)\|_2\leq C$, we can get
$$\|v(s)\|_{p_1+2}^{p_1+2}\leq C,\quad \|v(s)\|_{p_2+2}^{p_2+2}\leq C \quad {\rm for \ all }\quad  s\in [0,1).$$

In subcase (ii) $\lambda_1>0$, $\lambda_2\geq 0$, $\frac{2}{N}<p_1<p_2\leq \frac{4}{N}$.

We prove the conclusions in two steps.

Step 1. Let
\begin{align}
M(s)&=\frac{\lambda_1}{p_1+2}\int^s_0(1-\tau)^{\frac{Np_1-4}{2}-1}\int_{\mathbb{R}^N}|v(\tau)|^{p_1+2}dyd\tau,\label{110472}\\
N(s)&=\frac{\lambda_2}{p_2+2}\int^s_0(1-\tau)^{\frac{Np_2-4}{2}-1}\int_{\mathbb{R}^N}|v(\tau)|^{p_2+2}dyd\tau,\label{110473}\\
K(s)&=\frac{1}{2}\|\nabla v(s)\|_2^2+(1-s)[M'(s)+N'(s)].\label{1104101}
\end{align}
Then (\ref{10411}) can be written as
\begin{align}
K(s)=aM(s)+bN(s)+C_0\label{110481}
\end{align}
with
$$
a=\frac{4-Np_1}{2}, \quad b=\frac{4-Np_2}{2}
$$
and
$$
C_0=\frac{1}{2}\int_{\mathbb{R}^N}|\nabla \psi|^2dy+\frac{\lambda_1}{p_1+2}\int_{\mathbb{R}^N}|\psi|^{p_1+2}dy
+\frac{\lambda_2}{p_2+2}\int_{\mathbb{R}^N}|\psi|^{p_2+2}dy.
$$
Using (\ref{110481}), we have
$$(1-s)[aM'(s)+bN'(s)]\leq [aM(s)+bN(s)]+C_0.$$
Applying Gronwall's lemma, we have
\begin{align}
aM(s)+bN(s)\leq C(1-s)^{-\frac{4-Np_1}{2}}.\label{114133}
\end{align}

Step 2. Since $\frac{2}{N}<p_1<\frac{4}{N}$, we can take a constant $\epsilon_0$ satisfying
\begin{align}
0<\frac{4-Np_1}{2}<\epsilon_0<1\label{114131}
\end{align}
and we have
\begin{align}
K(s)+\int_0^s\frac{aM(\tau)+bN(\tau)+C_0}{(1-\tau)^{\epsilon_0}}d\tau
=aM(s)+bN(s)+C_0+\int_0^s\frac{K(\tau)}{(1-\tau)^{\epsilon_0}}.\label{1104102}
\end{align}
From (\ref{110481}), we know that $M''(s)$ and $N''(s)$ can be defined and are continuous in $[0,1)$. Consequently, $\left(\frac{aM(s)+bN(s)}{(1-s)^{\epsilon_0}}\right)''$ is also continuous in $[0,1)$. Hence $$\lim_{s\rightarrow 1}[\left(\frac{aM(s)+bN(s)}{(1-s)^{\epsilon_0}}\right)''-\frac{aM(s)+bN(s)}{(1-s)^{\epsilon_0}}]$$ exists(maybe equal to $+\infty$ or $-\infty$). First we prove that $$\lim_{s\rightarrow 1}[\left(\frac{aM(s)+bN(s)}{(1-s)^{\epsilon_0}}\right)''-\frac{aM(s)+bN(s)}{(1-s)^{\epsilon_0}}]\leq 0.$$  If $$\lim_{s\rightarrow 1}[\left(\frac{aM(s)+bN(s)}{(1-s)^{\epsilon_0}}\right)''-\frac{aM(s)+bN(s)}{(1-s)^{\epsilon_0}}]>0,$$ then there exists a $s_1$ such that
\begin{align}
\left(\frac{aM(s)+bN(s)}{(1-s)^{\epsilon_0}}\right)''>\frac{aM(s)+bN(s)}{(1-s)^{\epsilon_0}} \quad {\rm for} \quad s_1\leq s<1.\label{11411x1}
\end{align}
Using (\ref{11411x1}) and noticing that $M'(s)>0$ and $N'(s)>0$, we obtain
$$
\frac{aM(s)+bN(s)}{(1-s)^{\epsilon_0}}\geq C[e^s-e^{-s}]+C \quad {\rm for} \quad s_1\leq s<1.
$$
Consequently,
\begin{align}
aM(s)+bN(s)\geq C(1-s)^{-\epsilon_0}>>C(1-s)^{-\frac{4-Np_1}{2}} \quad {\rm for} \quad s_1\leq s<1,
\end{align}
which is a contradiction to (\ref{114133}). Hence $$\lim_{s\rightarrow 1}[\left(\frac{aM(s)+bN(s)}{(1-s)^{\epsilon_0}}\right)''-\frac{aM(s)+bN(s)}{(1-s)^{\epsilon_0}}]\leq 0.$$
And there exists $s_2\in [\frac{1}{2}, 1)$ such that
\begin{align}
\left(\frac{aM(s)+bN(s)}{(1-s)^{\epsilon_0}}\right)''\leq \frac{aM(s)+bN(s)}{(1-s)^{\epsilon_0}}
\quad {\rm for} \quad s_2\leq s<1.\label{11413x1}
\end{align}

On the other hand, using (\ref{11413x1}), after some elementary computations, we obtain
\begin{align}
&\qquad\int_0^s\frac{aM(\tau)+bN(\tau)+C_0}{(1-\tau)^{\epsilon_0}}d\tau\nonumber\\
&=-\frac{(1-s)^{1-\epsilon_0}[aM(s)+bN(s)+C_0]}{1-\epsilon_0}+\frac{C_0}{1-\epsilon_0}
+\int_0^s\frac{(1-\tau)^{1-\epsilon_0}[aM'(\tau)+bN'(\tau)]}{1-\epsilon_0}d\tau\nonumber\\
&=-\frac{(1-s)^{1-\epsilon_0}[aM(s)+bN(s)+C_0]}{1-\epsilon_0}+\frac{C_0}{1-\epsilon_0}
+\frac{s(1-s)^{1-\epsilon_0}[aM'(s)+bN'(s)]}{1-\epsilon_0}\nonumber\\
&\quad-\int_0^s\frac{\tau}{(1-\epsilon_0)}\left\{(1-\tau)^{1-\epsilon_0}[aM'(\tau)+bN'(\tau)]
\right\}'d\tau\nonumber\\
&=-\frac{(1-s)^{1-\epsilon_0}[aM(s)+bN(s)+C_0]}{1-\epsilon_0}+\frac{C_0}{1-\epsilon_0}
+\frac{s(1-s)^{1-\epsilon_0}[aM'(s)+bN'(s)]}{1-\epsilon_0}\nonumber\\
&\quad+\frac{\int_0^s\tau\left((1-\epsilon_0)(1-\tau)^{-\epsilon_0}[aM'(\tau)+bN'(\tau)]
-(1-\tau)^{1-\epsilon_0}[aM''(s)+N''(s)]\right)d\tau}{1-\epsilon_0}\nonumber\\
&=-\frac{(1-s)^{1-\epsilon_0}[aM(s)+bN(s)+C_0]}{1-\epsilon_0}
+\frac{C_0}{1-\epsilon_0}+\frac{s(1-s)^{1-\epsilon_0}[aM'(s)+bN'(s)]}{1-\epsilon_0}\nonumber\\
&\quad-\int_0^s\frac{\tau(1-\tau)}{1-\epsilon_0}\left(\frac{aM(\tau)+bN(\tau)}{(1-\tau)^{\epsilon_0}}\right)''d\tau
+\int_0^s\tau\left(\frac{aM(\tau)+bN(\tau)}{(1-\tau)^{\epsilon_0}}\right)'d\tau\nonumber\\
&\quad+\int_0^s\frac{2\epsilon_0\tau[aM'(\tau)+bN'(\tau)]}{(1-\epsilon_0)(1-\tau)^{\epsilon_0}}d\tau
+\int_0^s\frac{2\epsilon_0^2\tau[aM(\tau)+bN(\tau)]}{(1-\epsilon_0)(1-\tau)^{\epsilon_0+1}}d\tau\nonumber\\
&\geq -\frac{(1-s)^{1-\epsilon_0}[aM(s)+bN(s)+C_0]}{1-\epsilon_0}
+\frac{C_0}{1-\epsilon_0}+\frac{s(1-s)^{1-\epsilon_0}[aM'(s)+bN'(s)]}{1-\epsilon_0}\nonumber\\
&\quad-\int_0^s\frac{\tau(1-\tau)}{1-\epsilon_0}\times\frac{[aM(\tau)+bN(\tau)]}{(1-\tau)^{\epsilon_0}}d\tau
+\int_0^s\tau\left(\frac{aM(\tau)+bN(\tau)}{(1-\tau)^{\epsilon_0}}\right)'d\tau\nonumber\\
&\quad+\int_0^s\frac{2\epsilon_0\tau[aM'(\tau)+bN'(\tau)]}{(1-\epsilon_0)(1-\tau)^{\epsilon_0}}d\tau
+\int_0^s\frac{2\epsilon_0^2\tau[aM(\tau)+bN(\tau)]}{(1-\epsilon_0)(1-\tau)^{\epsilon_0+1}}d\tau\nonumber\\
&\geq -\frac{(1-s)^{1-\epsilon_0}[aM(s)+bN(s)+C_0]}{1-\epsilon_0}
+\frac{s(1-s)^{1-\epsilon_0}[aM'(s)+bN'(s)]}{1-\epsilon_0}\nonumber\\
&\quad-\int_0^s\frac{\tau(1-\tau)^{1-\epsilon_0}[aM(\tau)+bN(\tau)]}{1-\epsilon_0}d\tau
+C\quad {\rm for }\quad s_2<s<1.\label{1104103}
\end{align}
Using (\ref{114133})--(\ref{1104102}) and (\ref{1104103}), we get
\begin{align}
&\quad K(s)
+(1-s)^{1-\epsilon_0}[aM'(s)+bN'(s)]\nonumber\\
&\leq 2\int_0^s\frac{K(\tau)}{(1-\tau)^{\epsilon_0}}+2[aM(s)+bN(s)]
+\int_0^s\frac{(1-\tau)^{1-\epsilon_0}[aM(\tau)+bN(\tau)]}{(1-\epsilon_0)}d\tau+C\nonumber\\
&\leq 2\int_0^s\frac{K(\tau)}{(1-\tau)^{\epsilon_0}}+2[aM(s)+bN(s)]+C\quad {\rm for }\quad s_2<s<1.\label{1104112}
\end{align}
Taking $\kappa=\max(\epsilon_0, 1-\epsilon_0)$ and using (\ref{1104112}), we have
\begin{align}
(1-s)^{\kappa}J'(s)\leq 2J(s)+C\quad {\rm for }\quad s_2<s<1.\label{1104113}
\end{align}
with
\begin{align*}
J(s)=2\int_0^s\frac{K(\tau)}{(1-\tau)^{\epsilon_0}}+2[aM(s)+bN(s)].
\end{align*}
Applying Gronwall's lemma to (\ref{1104113}), we can obtain
\begin{align}
J(s)\leq C.\label{1104114}
\end{align}
Consequently,
\begin{align}
\|\nabla v(s)\|_2^2\leq 2J(s)+C\leq C'.\label{11411x5}
\end{align}
Using (\ref{11411x5}) and Gagliardo-Nirenberg's inequality, we have
$$
\|\nabla v(s)\|_{L^2}^2\leq C,\quad \|v(s)\|_{L^{p_1+2}}^{p_1+2}\leq C,\quad
\|v(s)\|_{L^{p_2+2}}^{p_2+2}\leq C $$
for all $s\in [0,1)$.

In subcase (iii) $\lambda_1>0$, $\lambda_2\geq 0$, $\frac{2}{N}<p_1\leq \frac{4}{N}<p_2\leq \frac{4}{N-2}$,
by ({\ref{110481}), we have
$$
\frac{1}{2}\|\nabla v(s)\|_2^2+(1-s)M'(s)+(1-s)N'(s)+|b|N(s)=aM(s)+C_0.
$$
Using the same technique as that in subcase (ii), we can get the conclusions of the proposition. We omit the details here.

Case (B) $\lambda_1>0$, $\lambda_2<0$, $\frac{2}{N}<p_1<p_2<\frac{4}{N}$.

If (\ref{104104}) is true,  we can take a $\varepsilon$ such that
$0<\varepsilon<\frac{Np_1-2}{2}$ and
$$\|\varphi(x)\|_{L^2}^{\frac{4}{N}}
<\frac{(4-Np_1)}{2N(p_2-p_1)C_N}\left(\frac{p_2+2}{|\lambda_2|}\right)^{\frac{4-Np_1}{N(p_2-p_1)}}
\left(\frac{\varepsilon\lambda_1(4-Np_1)}{(4-Np_2)(p_1+2)}\right)^{\frac{4-Np_2}{N(p_2-p_1)}}.$$
From (\ref{10411}) and using Young's inequality, we can get
\begin{align}
&\qquad\frac{1}{2}\int_{\mathbb{R}^N}|\nabla v(s)|^2dy+\frac{\lambda_1(1-s)^{\frac{Np_1-4}{2}}}{p_1+2}\int_{\mathbb{R}^N}|v(s)|^{p_1+2}dy\nonumber\\
&\leq
\frac{|\lambda_2|(1-s)^{\frac{Np_2-4}{2}}}{p_2+2}\int_{\mathbb{R}^N}|v(s)|^{p_2+2}dy
+C\nonumber\\
&\qquad -\frac{\lambda_1(Np_1-4)}{2(p_1+2)}\int^s_0(1-\tau)^{\frac{Np_1-4}{2}-1}\int_{\mathbb{R}^N}|v(\tau)|^{p_1+2}dyd\tau\nonumber\\
&\leq\frac{N(p_2-p_1)}{(4-Np_1)}\left(\frac{|\lambda_2|}{p_2+2}\right)^{\frac{4-Np_1}{N(p_2-p_1)}}
\left(\frac{(4-Np_2)(p_1+2)}{\varepsilon\lambda_1(4-Np_1)}\right)^{\frac{4-Np_2}{N(p_2-p_1)}}
\int_{\mathbb{R}^N}|v(s)|^{\frac{4}{N}+2}dy\nonumber\\
&\qquad+
\frac{\varepsilon\lambda_1(1-s)^{\frac{Np_1-4}{2}}}{(p_1+2)}\int_{\mathbb{R}^N}|v(s)|^{p_1+2}dy+C\nonumber\\
&\qquad -\frac{\lambda_1(Np_1-4)}{2(p_1+2)}\int^s_0(1-\tau)^{\frac{Np_1-4}{2}-1}\int_{\mathbb{R}^N}|v(\tau)|^{p_1+2}dyd\tau.
\label{10461}\end{align}
Using Gagliardo-Nirenberg's inequality, we have
\begin{align}
\int_{\mathbb{R}^N}|v(s)|^{\frac{4}{N}+2}dy&\leq C_N
\left(\int_{\mathbb{R}^N}|\nabla v(s)|^2dy\right)
\left(\int_{\mathbb{R}^N}|v(s)|^2dy\right)^{\frac{2}{N}}.\label{10462}
\end{align}
Since $\|v(s)\|_{L^2}=\|\psi(y)\|_{L^2}=\|\varphi(x)\|_{L^2}$,
then (\ref{10461}) and (\ref{10462}) imply that
\begin{align}
&\qquad c\int_{\mathbb{R}^N}|\nabla v(s)|^2dy+\frac{(1-\varepsilon)\lambda_1(1-s)^{\frac{Np_1-4}{2}}}{(p_1+2)}\int_{\mathbb{R}^N}|v(s)|^{p_1+2}dy
\nonumber\\
&\leq  \frac{\lambda_1(4-Np_1)}{2(p_1+2)}\int^s_0(1-\tau)^{\frac{Np_1-4}{2}-1}\int_{\mathbb{R}^N}|v(\tau)|^{p_1+2}dyd\tau+C.\label{10463}
\end{align}
Letting
$$\eta(s)=\frac{\lambda_1}{(p_1+2)}
\int^s_0(1-\tau)^{\frac{Np_1-4}{2}-1}\int_{\mathbb{R}^N}|v(\tau)|^{p_1+2}dyd\tau.$$
From (\ref{10463}), we have
\begin{align}
c\int_{\mathbb{R}^N}|\nabla v(s)|^2dy+(1-s)\eta'(s)\leq \frac{(4-Np_1)}{2(1-\varepsilon)}\eta(s)+C.\label{11410w1}
\end{align}
Using the comparison principle of ODE, and the relationship between the solution of in (\ref{11410w1}) and that of equation $$c\int_{\mathbb{R}^N}|\nabla v(s)|^2dy+(1-s)\eta'(s)=\frac{(4-Np_1)}{2(1-\varepsilon)}\eta(s)+C,$$
and noticing $\frac{4-Np_1}{2(1-\varepsilon)}<1$, similar to the arguments in case (A), we can obtain
\begin{align}
\|\nabla v(s)\|_{L^2}^2\leq C, \quad \|v(s)\|_{L^{p_1+2}}^{p_1+2}\leq  C,
\quad \|v(s\|_{L^{p_2+2}}^{p_2+2}\leq C\label{11223}\end{align}
for all $s\in [0,1)$. We omit the details here.\hfill $\Box$

By the results of Proposition 4.4, $\|v(s)\|_{H^1}$ is bounded as $s\rightarrow 1$. By Proposition 2.2, similar to the standard arguments as Theorem 7.4.1, Theorem 7.5.10 and Theorem 7.5.11 in \cite{Cazenave2}, we can establish the scattering theory in $\Sigma$ for (\ref{system0}) with $\lambda_1>0$. We omit the standard details here.

{\bf The proof of Theorem 2:} By the results of Subsection 4.1, the conclusions of Theorem 2 are true in cases (1) and (2). By the arguments in Subsection 4.2, the conclusions of Theorem 2 are true in cases (3) and (4). \hfill $\Box$

{\bf The proof of Corollary  1.1:} The corollary is a direct consequence of Theorem 2 in case (3) with $\lambda_2=0$.\hfill $\Box$

{\bf The proof of Theorem 3:} In the proof of Proposition 4.4, we obtain
$$
\|\nabla v(s)\|_2^2\leq C \quad{\rm for \ all } \quad s\in [0,1).
$$
if (i) $\lambda_1>0$, $\lambda_2\geq 0$, $\frac{2}{N}<p_1<p_2\leq \frac{4}{N-2}$ or (ii) $\lambda_1>0$, $\lambda_2<0$, $\frac{2}{N}<p_1<p_2<\frac{4}{N}$ with (\ref{104104}). Hence for any $2\leq r\leq \frac{2N}{N-2}$, we have
$$\|v(s)\|_{L^r}\leq C.$$
 Noticing that $(1+t)^{\frac{N\beta}{2}}\|u(t)\|_{\beta+2}^{\beta+2}=\|v(s)\|_{\beta+2}^{\beta+2}$ for any $\beta\geq 0$, we obtain
some results on the decay of solutions to (\ref{system0}). That is,
\begin{align}
\|u(t)\|_{L^r}\leq C(1+t)^{-\frac{N(r-2)}{2r}}\label{113121}
\end{align}
for all $t\in [0, +\infty)$.\hfill $\Box$

\section{Comments}
\qquad In the last section, we want to give some comments on our results.

{\bf Comment 5.1.} Although the techniques are elementary in this paper, we give some results on the scattering for (\ref{system0}) in case (I) and (II), which partly solves some open problems on scattering for (\ref{system0}).
By the way, as a direct consequence
of our results, we establish the scattering for (\ref{classic1}) with $\lambda>0$ and $\frac{2}{N}<p<\alpha_0$, which is also an open problem in this direction.

We also would like to compare our results with those of \cite{Tao}.

(1) We establish the no scattering results for (\ref{system0}) with $p_1\leq \frac{2}{N}$, which is not covered by those in \cite{Tao}.

(2) Our results contain the scattering for (\ref{system0}) in case (i) $\lambda_1>0$, $\lambda_2\geq 0$, $\frac{2}{N}<p_1\leq\alpha_0$ and $0<p_2\leq \frac{4}{N-2}$ and case (ii) $\lambda_1>0$, $\lambda_2<0$, $\frac{2}{N}<p_1<p_2<\frac{4}{N}$, which is not covered by those in \cite{Tao}.

(3) Our results contain the scattering for (\ref{system0}) in case (i) $\lambda_1<0$, $\lambda_2\geq 0$, $\frac{4}{N+2}<p_1<p_2<\frac{4}{N-2}$ and case (ii) $\lambda_1<0$, $\lambda_2<0$, $\frac{4}{N}<p_1<p_2<\frac{4}{N}$, which is not covered by those in \cite{Tao}.

(4) The scattering for (\ref{system0}) with
$\lambda_1>0$, $\lambda_2\geq 0$, $\alpha_0<p_1<p_2\leq \frac{4}{N-2}$ is established both in this paper and \cite{Tao}. However, our methods rely on the elementary technique(Proposition 4.4), while theirs rely on traditional Strichartz's estimate for $(x+2it\nabla)u$.

(5) Tao, Visan and Zhang had established the scattering theory in $H^1$ for (\ref{system0}) in \cite{Tao}, however, we cannot obtain any results on the scattering theory in $H^1$ for (\ref{system0}) in this paper. We think that every method has its weakness. We also cannot use the method here to get
more scattering properties for (\ref{system0}) with  $\lambda_1<0$. The scattering theory hasn't been established in case (i) $\lambda_1<0$, $\lambda_2>0$, $p_1>\frac{4}{N+2}$ with $p_2=\frac{4}{N-2}$ or case (ii) $\lambda_1<0$, $\lambda_2>0$, $0<p_1\leq \frac{4}{N+2}$ with $0<p_2\leq \frac{4}{N-2}$. However, we suspect that there exists scattering sate for (\ref{system0}) at least in in $L^2$ topology in case (i), while there are no scattering theory for (\ref{system0}) even in in $L^2$ topology in case (ii).

{\bf Comment 5.2.} We would like to compare Theorem 1 with Theorem 7.5.2 and Remark 7.5.5 in \cite{Cazenave2}.
 The results of Theorem 1 show that:  If $p_1<\frac{2}{N}$, the role of $-\mu|u|^{p_1}u$ prevails that of $\lambda_2|u|^{p_2}u$, both Theorem 7.5.4 in \cite{Cazenave2} and Theorem 1 illustrate that the power $\frac{2}{N}$ can be look as a border of wether the solution has scattering state or not. On the other hand, since the interaction between the defocusing nonlinearity and the focusing one, our conclusions are differ from those of Remark 7.5.5 (ii) in \cite{Cazenave2}(Some scattering results on (\ref{classic1}) with the focusing nonlinearity were given there).

{\bf Comment 5.3.}  We would like to compare Theorem 2 with those theorems in Section 7.5 of \cite{Cazenave2}. Theorem 2 shows that: If $\frac{2}{N}<p_1$, the role of nonlinearity $\lambda_1|u|^{p_1}u$ overwhelm that of $\lambda_2|u|^{p_2}u$ under some assumptions, and we can look the nonlinearity $\lambda_2|u|^{p_2}u$ as a disturbance. Therefore, if $\lambda_1>0$, $p_1>\alpha_0$ or $\lambda_1<0$, $p_1>\frac{4}{N+2}$, our results meet with those of \cite{Cazenave2}. However, we give the information on the scattering for (\ref{classic1}) with $\lambda>0$, $\frac{2}{N}<p_1\leq \alpha_0$.

{\bf Comment 5.4.} By the results of this paper and \cite{Tao}, the nonlinearity $\lambda_1|u|^{p_1}u$ has the main role of scattering for (\ref{system0}). On the other hand, the focusing nonlinearity $\lambda_2|u|^{p_2}u$(i.e.,  $\lambda_2<0$) may lead the phenomenon of finite time blowup for (\ref{system0}) happen.

{\bf Comment 5.5.} The method in this paper and those of \cite{Tao} can be used to deal with the following Cauchy problem
\begin{align}
\label{gs} \left\{\begin{array}{ll}&i u_t+\Delta
u=\sum_{i=1}^m\lambda_i|u|^{p_i}u, \quad x\in\mathbb{R}^N,
\quad t\in\mathbb{R}\cr
 &u(0,x)=\varphi(x), \quad x\in\mathbb{R}^N,
\end{array}
\right. \end{align} where $N\geq 3$, $0<p_1<p_2<...<p_m\leq\frac{4}{N-2}$,
$\lambda_i$, $i=1,2,...,m$  are real constants. In many cases, whether the solution of (\ref{gs}) possess a scattering state or not are essentially depended on the nonlinearities $\lambda_1|u|^{p_1}u$ and $\lambda_m|u|^{p_m}u$, because $\lambda_i|u|^{p_i}u$, $i=2,...,(m-1)$ can be controlled by $\lambda_1|u|^{p_1}u$ and $\lambda_m|u|^{p_m}u$ if one use Young's inequality.
~~~~~~~~~~~~~~~~~~~~~~\\
~~~~~~~~~~~~~~~~~~~~\\
{\bf Acknowledgement}

This work is supported by the National Nature Science Foundation
of China(Project  No. 11071237).

\end{document}